%% file: nicola.tex
\documentclass[graybox]{svmult}

\usepackage{type1cm}        % activate if the above 3 fonts are
\usepackage{makeidx}         % allows index generation
\usepackage{graphicx}        % standard LaTeX graphics tool
\usepackage{multicol}        % used for the two-column index
\usepackage[bottom]{footmisc}% places footnotes at page bottom

\usepackage{newtxtext}       % 
\usepackage{newtxmath}       % selects Times Roman as basic font

\makeindex             % used for the subject index
\begin{document}

\title*{A life in Mathematical Analysis: a conversation with Luigi Rodino}
\author{Fabio Nicola}
\institute{Fabio Nicola \at Dipartimento di Scienze Matematiche, Politecnico di Torino, Corso Duca degli Abruzzi 24, 10129 Torino, Italy, \email{fabio.nicola@polito.it}}
\maketitle

\abstract*{Each chapter should be preceded by an abstract (no more than 200 words) that summarizes the content. The abstract will appear \textit{online} at \url{www.SpringerLink.com} and be available with unrestricted access. This allows unregistered users to read the abstract as a teaser for the complete chapter.
Please use the 'starred' version of the \texttt{abstract} command for typesetting the text of the online abstracts (cf. source file of this chapter template \texttt{abstract}) and include them with the source files of your manuscript. Use the plain \texttt{abstract} command if the abstract is also to appear in the printed version of the book.}

\abstract{This note is the transcription of an interview with Professor Luigi Rodino, on the occasion of the ISAAC-ICMAM Conference of Analysis in Developing Countries (December 2, 2024 – Bogotà), that was dedicated to him. Luigi Rodino is at present Emeritus Professor at the University of Turin, and a member of the Accademia delle Scienze di Torino.}

\section*{}
{\bf Fabio:} Welcome Luigi, and thank you for accepting to do this interview. 

You were born in 1948 in Cairo Montenotte, a small town in Savona province. Your father was a physician, your mother worked as a teacher and liked writing. 

Has your passion for mathematics grown over time or does it date back to when you were a child?
\par\medskip
\noindent {\bf Luigi:}  
My father had in mind for me the studies and profession in medicine and while I was a child he was training me in this direction, explaining to me the scientific and human aspects of the different diseases. He was certainly disappointed in my choice of mathematics, though later he appreciated my academic career. Instead, my mother approved from the beginning my choice, because of the cultural aspects of the research and teaching activities.

The persons who addressed me to the choice were my professors in mathematics in high school, Urbani and Spotorno. Beside excellent teachers, they had rather non-standard temperaments. Mrs. Urbani, a strong woman, had a relevant role as a warrior in the communist partisans during the second world-war. Bruno Spotorno was simultaneously a high school teacher and professor of education problems in the University of Genova, he had eccentric and amusing behaviors (as people sometimes expect from a mathematician!).
Let me say frankly:  when entering the University I did not really have a passion for mathematics, but after my decision I intended to do my best in this direction.         
\par\medskip

\noindent{\bf Fabio:}  You graduated in Mathematics at the University of Turin in 1971. As you told me, at that time, in Turin, there were no experts in distribution theory, and you acquired the basics of that theory through a course given by two former students of Laurent Schwartz – François Treves and Louis Boutet de Monvel – which significantly shaped your future path. Then, in 1972, you moved to Lund to study with Hörmander, where you authored your first paper, {\em A class of pseudo-differential operators on the product of two manifolds and applications} \cite{rod1}. That paper was well received by the scientific community and can be considered an early contribution to the vector-valued pseudodifferential calculus, later carried on by Louis Boutet de Monvel, Alain Grigis, Bernard Helffer, Bert-Wolfgang Schulze, Richard Melrose, and others. 

Do you want to share any anecdote about that period? 

\par\medskip
\noindent {\bf Luigi:} 
Let me first say a few words about the University of Turin during my studies, 1967-1971. Mathematical analysis was dominated there by F. G. Tricomi. He was still in service, but actually his activity was limited to a seminar. The courses were delivered by his assistants, following word by word several books he wrote: calculus, advanced calculus, ordinary and partial differential equations, complex variables, special functions. He was strongly against the modern functional analysis, he wrote somewhere "the theory of the Schwartz distributions is repellent". Nevertheless, this classical background turned out to be very useful to me in the subsequent studies.

Ph.D. programmes did not exist at that moment in Italy, and were replaced by the grant "assegno di ricerca", that I got at the Institute of Mathematics of the Politecnico of Torino. The head of the Institute, Piero Buzano, kindly allowed me to spend the grant abroad. Namely, under the suggestion of another professor of the Institute, Giuseppe Geymonat, I asked Lars H\"ormander to accept me as a student in Lund, Sweden. Jointly with the summer courses of Treves and Boutet de Monvel that you mention, H\"ormander opened my mind to mathematical research. I remember my first meeting in Lund with H\"ormander. He was extremely happy that I had read his 1971 paper on Fourier Integral Operators \cite{hor1971}, and we discussed applications of his definition of wave front set. His recent students J. Sj\"ostrand and A. Melin were in France and Denmark at that moment, and other former students of H\"ormander and Gårding were working on constant coefficient operators. So, I was lucky in some sense to have special attention to me in the new directions. His help was invaluable.

I still keep the manuscript of my first papers, in particular the one you mention, with precise hand-written corrections of H\"ormander. I remember my stay in Sweden as a wonderful period, in particular, for the friendship of the people there.

\par\medskip

\noindent{\bf Fabio:} After a period spent at the University of Princeton with Joseph J. Kohn, you became a full professor at the University of Turin in 1981. Your research activity covered a number of topics, including the problems of hypoellipticity and local solvability, the theory of pseudodifferential operators in the Gevrey category, and in more recent years, global pseudodifferential calculus and time-frequency analysis. 

Were there any people, perhaps collaborators, who particularly influenced your research in the early stages of your career?

\par\medskip
\noindent {\bf Luigi:} 
After my stay in Sweden, I was for one year in Princeton with my wife. The scientific environment was excellent but, compared with our stay in Sweden, we had some difficulties in adapting to the US life, especially concerning food! Meanwhile I turned attention to the problems of hypoellipticity and local solvability. Under the suggestion of Boutet de Monvel I collaborated on these problems with Bernard Helffer and, when finally returning to Italy, with Cesare Parenti from Bologna University. In fact it was natural for me to join people in Italy working on H\"ormander's setting, as we now say "microlocal analysis". So I also joined the research group of Lamberto Cattabriga, studying Gevrey classes with pseudo-differential techniques.

\begin{figure}[h]
\sidecaption
% For example, with the graphicx style use
\includegraphics[scale=.55]{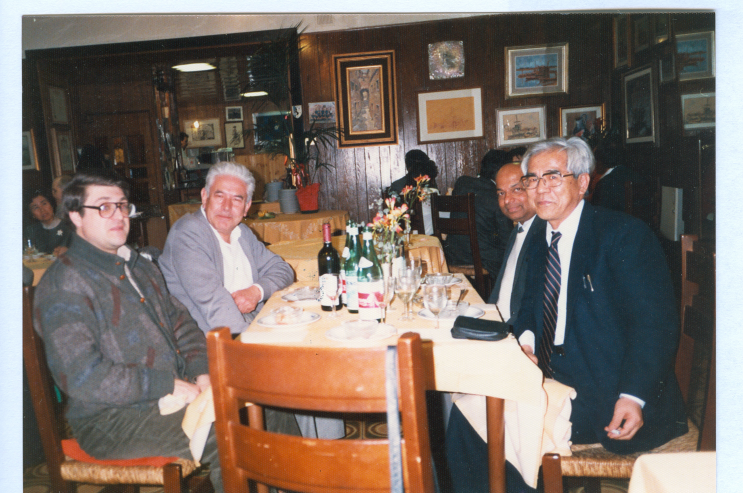}
%
% If no graphics program available, insert a blank space i.e. use
%\picplace{5cm}{2cm} % Give the correct figure height and width in cm
%
\caption{Luigi Rodino and Ennio De Giorgi (left), Venkatesha Murthy
e  Sigeru Mizohata (right)}
\label{figura_degiorgi}       % Give a unique labe
\end{figure}

Through Cattabriga I met his collaborator and close friend Ennio De Giorgi (Figure \ref{figura_degiorgi}). De Giorgi is now recognized as the best italian mathematician of the second-half of the past century. I remember him as an outstanding, limpid mind, kind sweet person and, in harmony with my feelings, a gourmet and an excellent connoisseur of wines.
It is hard to mention all my collaborators and all the people which had an influence on my career. Coming to recent times, I must say that young collaborators and all my Ph.D. students influenced me a lot. It seems strange, but young people are more sensitive to the new trends. They forced me into emerging directions of research. This is the case for the fields you mention, global pseudo-differential calculus and time-frequency analysis, which are in my present interest of research.
\par\medskip

\noindent{\bf Fabio:}  Over the years, you have taken on several administrative roles at the University of Turin, such as director of the department, president of the Mathematics Program, and Ph.D. Program coordinator. You coordinated several national and international research projects and were president of the ISAAC society. At present, you are a member of the Accademia delle Scienze di Torino and a member of the working group of the Italian Mathematical Union dedicated to mathematical collaboration projects with  the  Global South. Alongside Cesare Parenti, you are recognized as the founder of the Italian school of microlocal analysis. In the early 2000s, you also played a significant role in revitalizing this field by founding the Journal of Pseudo-Differential Operators and Applications, together with Man Wah Wong. 

Recently, I saw a renewed interest in microlocal analysis, especially in the mathematical physics community (quantum field theory on curved spaces, etc.), and also in the time-frequency analysis community.  How do you see the research landscape within mathematical analysis? Are there, in your opinion, particularly promising lines of research?

\par\medskip
\noindent {\bf Luigi:} 
In mathematics, as in other sciences, there are people who are fond of their field and obtain excellent results by an exclusive dedication, without taking care of practical problems around. On the opposite side, there are people with a relatively small personal scientific production, who are excellent organizers, and help their community assuming directive and administrative duties in their institutions. I tried to compromise between the two options. A crucial moment was when I was head of the Department of Mathematics in the University of Torino for three years. During this period I had to stop my research, and I realized that choosing this kind of activities I had to renounce a fair research production.

Nevertheless I enjoyed acting as president of ISAAC.  As for other international associations, you have the chance to act in a political direction. Despite ISAAC being born in a very rich nation thanks to the original ideas of Robert Gilbert Chancellor of Delaware University, all the ISAAC presidents acted in favor of disadvantaged countries, promoting friendship among all world nations. The task is not difficult for mathematics, because we have a tradition of international exchange.

To answer your final question: yes, there is a renewed interest in microlocal analysis. With respect to the definition of the wave front set of H\"ormander, people are more addressed to quantitative formulations of microlocalization. This leads to a new deep analysis of the uncertainty principle, and applications in signal analysis and quantum physics.

Concerning the general landscape within mathematical analysis, as editor of some journals I observed an extremely large production, in particular on applicative aspects. With respect to the preceding literature, coming mainly from the US and Europe, one notes the appearance of new nationalities with an impressive number of crossed collaborations. The enlargement of the role of mathematics in world society is certainly a positive fact.

\par\medskip

\noindent{\bf Fabio:}  Your academic contributions include supervising 17 Ph.D. students and authoring six books along with approximately 150 scientific papers, receiving about 2000 citations. 

Is there any mathematical problem that you would have liked to tackle, and that ended up in a drawer? You know, a bee in your bonnet, that you think about from time to time.

\par\medskip
\noindent {\bf Luigi:} 
In the past I had some unreasonable hopes concerning extremely difficult questions in the general theory of the linear partial differential operators, in particular for the hypoellipticity problem. In the paper entitled "Pseudo differential operators and non-elliptic boundary problems" 1966 \cite{hor1966}, H\"ormander gave general implicit necessary and sufficient conditions for hypoellipticity with fixed loss of Sobolev derivatives. More explicit results were given later by  F. Treves, V. V. Grushin, Y. V. Egorov and others. My collaboration with Helffer and Parenti was in this direction. I had a secret hope, to solve in general the hypoellipticity problem by testing the equation on classes of functions with singularities, without the use of Fourier analysis. This remained a dream.
A related, apparently easier problem, is the global hypoellipticity in the Schwartz spaces $S$ and $S'$ of the operators with polynomial coefficients. You know Fabio, we have a paper together concerning the one dimensional case. I had a general conjecture that I proposed in some international meetings, in Japan, Mexico, etc., around 10 years ago. The audiences were suspicious, and they were right, because the conjecture was wrong. 
At present, after retirement, I tackle only reasonable problems that I feel confident are solvable.
\par\medskip

\noindent{\bf Fabio:}  During my freshman year at university, I had the pleasure of meeting you while attending the Mathematical Analysis I class. Your teaching instantly captivated me because of two key characteristics: exemplary clarity and a constant hint of irony. Attending your lectures was indeed a pleasure. Subsequently, I also took your course on partial differential equations, which was in fact inspired by the above-mentioned course given by François Treves and Louis Boutet de Monvel. 

I guess teaching was a real pleasure for you too, right?

\par\medskip
\noindent {\bf Luigi:}  
Yes, teaching was a pleasure for me. Besides courses for graduation and Ph.D. students, I gave for about 40 years courses of background analysis. Despite the contents being elementary, I prepared myself carefully. I must confess that I prepared as well some jokes or amusing stories to relieve the students after the proof of difficult theorems. Recently, a young professor of probability in the department recalled that, when following the first lecture of calculus at the university, and expecting stressing arguments, was surprised by the teacher, spending a quarter of hour on the discussion of different qualities of chocolates in Torino. I do not remember such a lecture, but I recognize myself as this teacher.
\par\medskip

\begin{figure}[h]
\sidecaption
% For example, with the graphicx style use
\includegraphics[scale=.25]{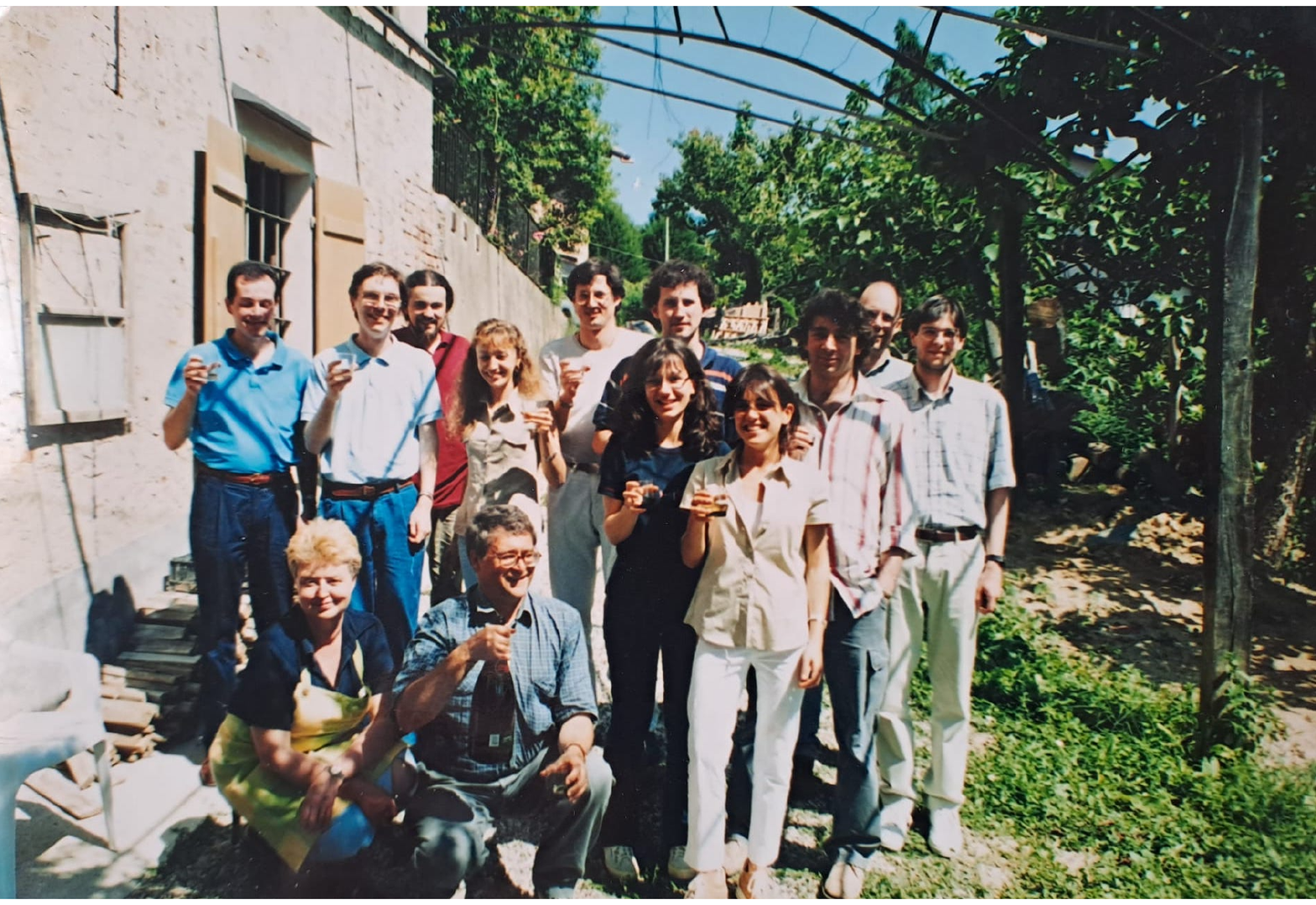}
%
% If no graphics program available, insert a blank space i.e. use
%\picplace{5cm}{2cm} % Give the correct figure height and width in cm
%
\caption{Luigi Rodino and his wife (center), with (former) students}
\label{figura_pranzo}       % Give a unique labe
\end{figure}

\noindent{\bf Fabio:}  Throughout our numerous scientific collaborations, I spent many hours in your office. At that time, you had numerous students and collaborators; entering your office felt akin to a doctor's office, and reservation was recommended. The community spirit in your research group was fostered by annual lunches you hosted at your country house, where you invited your students and close collaborators (Figure \ref{figura_pranzo}). I have also been particularly impressed by your constant concern for young researchers facing difficulties in getting a permanent position.

I guess the possibility of creating a school represented the greatest human and professional satisfaction of your career. What advice would you give to a young student who intends to become a professional mathematician?

\par\medskip
\noindent {\bf Luigi:}  
In Italy when graduating or getting a master's degree, after some years of university studies (three, four or five years, according to different rules in different periods), one is asked to produce an "honour thesis". This is not required to be a research work, though talented students may reach results worth publication. I accepted to supervise about 150 of these dissertations. Some colleagues told that I was crazy, accepting also students with very low marks, but I tell you, my regret is that I had to decline, because of other engagements, the applications of some people, having confidence on my help.

You refer to Ph.D. students. In Italy they are selected by a rather difficult examination, so you are sure to deal with bright people. Some of my Ph.D. Students are now professors of mathematics in different universities in Italy and abroad. Some of them are working in private enterprises or in the high school. All of them defended with success their dissertation, with full satisfaction from my part and their part.
Also nowadays I would encourage young students to an academic career, despite that the present opportunities are not so large. You say I am ironic, so let me give to them the following suggestions: try to publish immediately, also if your results are not complete; write as many papers as you can, possibly repeating the same result in different issues; pay extreme attention to impact factor of the journals. At the end let me be serious: enjoy mathematics, obtaining new results or new perspectives will give to you a lot of pleasure.
\par\medskip

\noindent{\bf Fabio:} Finally, I like to recall a kind of poster displayed on your office door, reading {\em Elogio della Mitezza} --- that is, {\em  Praise of Meekness} --- alongside a picture of a cat. Was meekness a choice for you, or did it always come naturally to you?

\par\medskip
\noindent {\bf Luigi:} 
Let me be extremely serious: we see in modern society some horrible differences between rich people and poor people, rich nations and poor nations. Even more horrifying, we see at this moment a series of wars all over the world, main motivation being the commerce of weapons of different types. Some well recognized politicians say that this will help economy. Disgusting!

As I get older, I feel more and more sad and angry about this. We must try our best to establish peace. I know you share my feelings. Thank you Fabio for this interview.
\par\medskip

\noindent{\bf Fabio:} Thank you again, Luigi, for this interview. Let me take this opportunity to express my gratitude, also on behalf of the other former students, for being an exceptional master and mentor, for your generosity and friendship.

\begin{acknowledgement}
Fabio Nicola and Luigi Rodino express their gratitude to the organizers of the ISAAC-ICMAM Conference of Analysis in Developing Countries (December 2, 2024 – Bogotà). They are also very grateful to Patrik Wahlberg for carefully reading a preliminary version of this note and suggesting several stylistic improvements. 

Fabio Nicola is a Fellow of the {\it Accademia delle Scienze di Torino} and a member of the {\it Societ\`a Italiana di Scienze e Tecnologie Quantistiche (SISTEQ)}. 
\end{acknowledgement}

\begin{figure}[h]
\sidecaption
% For example, with the graphicx style use
\includegraphics[scale=.25]{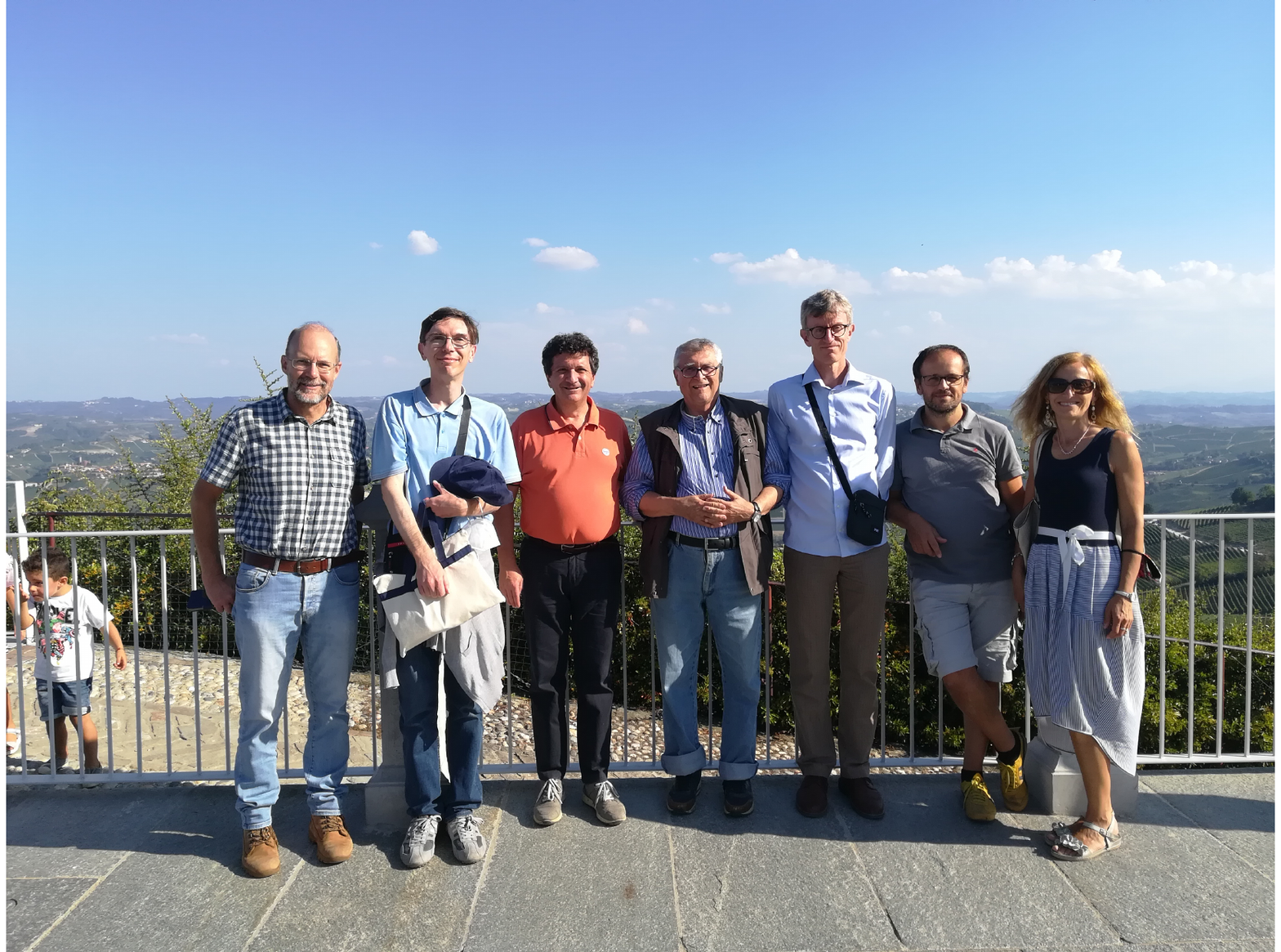}
%
% If no graphics program available, insert a blank space i.e. use
%\picplace{5cm}{2cm} % Give the correct figure height and width in cm
%
\caption{Luigi Rodino (center), with former students and collaborators}
\label{figura_convegno}       % Give a unique labe
\end{figure}

\vskip5 cm

\input{references}

\end{document}

%% file: references.tex
%%%%%%%%%%%%%%%%%%%%%%%% referenc.tex %%%%%%%%%%%%%%%%%%%%%%%%%%%%%%
% sample references
% %
% Use this file as a template for your own input.
%
%%%%%%%%%%%%%%%%%%%%%%%% Springer-Verlag %%%%%%%%%%%%%%%%%%%%%%%%%%
%
% Journal article